\documentclass[11pt,a4paper,reqno]{amsart} 

\usepackage{amssymb,graphicx}
      
\usepackage{amssymb}

\newcommand{\koniec}{\begin{flushright}  $\Box $ \end{flushright}}
\newtheorem{theo}{Theorem}[section] 
\newtheorem{prop}[theo]{Proposition}  
\newtheorem{lemma}[theo]{Lemma}
\newtheorem{defi}[theo]{Definition}

\newcounter{mnotecount}[section]

\renewcommand{\themnotecount}{\thesection.\arabic{mnotecount}}

\newcommand{\mnote}[1]
{\protect{\stepcounter{mnotecount}}$^{\mbox{\footnotesize
$
\bullet$\themnotecount}}$ \marginpar{
\raggedright\tiny\em
$\!\!\!\!\!\!\,\bullet$\themnotecount: #1} }

\newcommand{\hook}{{\setlength{\unitlength}{11pt}   
                   \begin{picture}(.833,.8)
                   \put(.15,.08){\line(1,0){.35}}
                   \put(.5,.08){\line(0,1){.5}}
                   \end{picture}}}
\newcommand{\CP}{\mathbb{CP}}
\newcommand{\C}{\mathbb{C}}
\newcommand{\PP}{\mathbb{P}}
\newcommand{\RP}{\mathbb{RP}}
\newcommand{\R}{\mathbb{R}}

\def\p{\partial}
\def\be{\begin{equation}}
\def\ee{\end{equation}}

\def\bea{\begin{eqnarray}}
\def\eea{\end{eqnarray}}

\newcommand{\spp}{\mathbb{S}}

\begin{document}\date{November 4, 2011}
\title{Co--calibrated $G_2$ structure from cuspidal cubics}
\author{Boris Doubrov}
\address{Belarussian State University\\ Nezavisimosti Av. 4, \\
220030, Minsk, Belarus}
\email{doubrov@islc.org}
\author{Maciej Dunajski}
\address{Department of Applied Mathematics and Theoretical Physics\\ 
University of Cambridge\\ Wilberforce Road, Cambridge CB3 0WA, UK.
}
\email{m.dunajski@damtp.cam.ac.uk}
\begin{abstract}
We establish a twistor correspondence between a cuspidal cubic curve 
in a complex projective plane, and a co-calibrated homogeneous $G_2$ structure on the seven--dimensional parameter space of such cubics. Imposing the Riemannian reality conditions leads to an explicit co-calibrated $G_2$ structure on $SU(2, 1)/U(1)$. This is an example of an $SO(3)$ structure in seven dimensions.

Cuspidal cubics and their higher degree analogues with constant projective curvature are characterised as integral curves of 7th order ODEs. Projective orbits of such curves are shown to be analytic continuations of Aloff--Wallach manifolds, and it is shown that
only cubics lift to a complete family of contact rational curves in a projectivised cotangent bundle to a projective plane.
\end{abstract}   
\maketitle
\section{Introduction}
Twistor theory gives rise to correspondences between global
algebraic geometry of rational curves in complex two--folds
or three--folds, and local differential geometry on the moduli spaces
of such curves. The embedding of a rational curve $L$ in a complex manifold
$\mathcal{Z}$ is, to the first order, described by the normal bundle
$N(L):=T\mathcal{Z}/TL$. This is a holomorphic vector bundle, and thus
by the Birkhoff--Grothedieck theorem it is a direct sum of
$(dim(\mathcal{Z})-1)$ line bundles $\mathcal{O}(n)$ of degree $n$
which can vary between the summands. The Kodaira deformation theorem
\cite{kodaira}  states that if $H^1(L, N(L))=0$, then $L$ belongs
to a locally complete family $\{L_m, m\in M\}$ where $M$ is some complex 
manifold, and there exists a canonical isomorphism
\[
T_m M\cong H^0(L_m, N(L_m)).
\]
In the original Non--linear Graviton construction
of Penrose \cite{Pe76} the twistor space $\mathcal{Z}$ is a complex 
three--fold and the normal bundle is $N(L)=\mathcal{O}(1)\oplus\mathcal{O}(1)$.
This gives rise to an anti--self--dual conformal structure on a 
four--dimensional manifold $M$ such that the null vectors fields in $M$
correspond to sections of $N(L)$ vanishing at one point. 
In the subsequent twistor constructions of Hitchin \cite{hitchin},
the twistor space $\mathcal{Z}$ is a two--fold, and $N(L)=\mathcal{O}(1)$
or $\mathcal{O}(2)$. In the first case $M$ is a surface admitting a projective structure, and in the second case $M$ is a three--dimensional manifold
with an Einstein--Weyl structure. 
  The whole set up can be 
generalised to contact rational curves in complex three--folds \cite{Bryant2}.
The moduli space of such curves with normal bundle 
$\mathcal{O}(n)\oplus\mathcal{O}(n)$ admits an integrable  
$GL(2)$ structure \cite{Bryant2, DT}. See \cite{Dbook} for other examples of twistor constructions.

 The aim of this paper is to use a twistor correspondence to construct seven
dimensional manifolds with $G_2$ structure. The general theory was
developed in \cite{DG10}, and in the present paper we construct a class
of explicit new examples corresponding to $L$ being a plane cuspidal cubic
\be
\label{cuspidal}
y^2-x^3=0
\ee
in a complex two--fold $\mathcal{Z}=\CP^2$. In order to do that, we need
to refine the twistor correspondence as outlined above, because
the cuspidal cubics, although rational, are singular in the complex 
projective plane.  This can be dealt with either 
by considering the contact lifts of the cuspidal cubics to 
$\PP(T\CP^2)$, where they become smooth contact curves with normal bundle
$\mathcal{O}(5)\oplus\mathcal{O}(5)$, or by working directly with singular curves. Both approaches lead to deformation theory
of (\ref{cuspidal}) as a cuspidal cubic curve (rather than as a general plane cubic). We shall find that  normal vector field 
to a cuspidal cubic $L$ vanishes at six general points on $L$ away from the cusp. The parameter space of cuspidal  cubics arising from this deformation theory is the seven-dimensional homogeneous space $M=PSL(3, \C)/\C^*$. 
All cuspidal cubics are projectively equivalent and belong to the same 
$PSL(3, \C)$ orbit of (\ref{cuspidal}) in $\CP^2$.

To formulate our main result, recall that the $(n+1)$ dimensional space
of holomorphic sections $H^0(\CP^1, \mathcal{O}(n))$ is isomorphic to
the vector space $\mathcal{V}_n=\mbox{Sym}^n(\C^2)$ of homogeneous polynomials of degree $n$ 
in two variables $(s, t)$. Any such section is of the form
\be
\label{sextic111}
V(s, t)=v_0t^{n}+nv_1 t^{n-1}s+\frac{1}{2}n(n-1) v_2 t^{n-2}s^2+
\dots +v_{n} s^{n}.
\ee
Let $U, V$ be elements
of $\mathcal{V}_n$. The $p$th transvectant of two polynomials
$U$ and $V$  is an element of  $\mathcal{V}_{2n-2k}$ given by
\[
<U, V>_p=\frac{1}{p!}\sum_{i=0}^p(-1)^i{p\choose i}
\frac{\p ^p U}{\p t^{p-i}\p s^i}\frac{\p ^p V}{\p t^{i}\p s^{p-i}}.
\]
We shall first establish (Proposition \ref{prop_GL}) a canonical
identification between vector fields on $M$ and elements of
$\mbox{Sym}^6(\C^2)$, and then use it to prove
\begin{theo}
\label{main_theo}
The seven--dimensional space $M=SL(3, \C)/\C^*$ of plane cuspidal
cubics admits a canonical complexified $G_2$ structure where the
three form $\phi$ and the metric $g$ 
\[
\phi(U, V, W)=<<U, V>_3,W>_3, \quad g(U, U)=<U, U>_6
\]
are explicitly given by {\em(}\ref{formula1}{\em)}. This $G_2$ structure is co--calibrated
i. e.
\be
\label{cocalibrated}
d\phi=\lambda *\phi+*\tau, \quad d*\phi=0,
\ee
where $\lambda$ is a constant and $\tau$ is a certain three--form
such that $\phi\wedge\tau=\phi\wedge *\tau=0$.
\end{theo}
We are ultimately interested in real $G_2$, and
we shall show that the structure from Theorem \ref{main_theo} admits 
three homogeneous  real forms: two with signature $(4, 3)$, 
where $M=SL(3, \R)/\R^*$ or  $M=SU(3)/U(1)$, and one Riemannian
with $M=SU(2, 1)/U(1)$.  

The paper is organised as follows: In the next Section we shall summarise
basic facts about the nonlinear group actions on spaces of symmetric polynomials, with the particular emphasis on sextics in two variables
and cubics in three variables. In Section \ref{section_gl2}
we shall demonstrate (Proposition \ref{prop_GL}) that the family $M$ of cuspidal cubics  admits a $GL(2)$ structure, so that the vectors tangent to $M$
can be identified with elements of $\mbox{Sym}^6(\C^2)$.
In Section \ref{section_main} we shall prove Theorem \ref{main_theo}
and construct  the conformal structure and the associated three-form 
directly from the $GL(2)$ structure on $M$.
Restricting to a real Riemannian slice will reveal
a co-calibrated $G_2$ structure on $SU(2, 1)/U(1)$. 
In Section \ref{section_de1}, we shall discuss the differential equations approach, where $M$ 
arises as the solution space  of a 7th order ODE. To find this ODE,
take seven derivatives of the general form of the cuspidal cubic and express the seven parameters  in terms of $y(x)$ and its first six derivatives, which leaves one condition in the form of 
the ODE (Proposition \ref{cubic_ode_prop}). From the point of view of projective geometry of curves,
the cuspidal cubics belong to the class of algebraic curves with constant projective curvature \cite{Wilczynski_book}.
We shall find (Proposition \ref{flat_curves}) all curves in $\CP^2$  
which give rise to  seven--dimensional projective orbits. 
These
curves have constant projective curvature and are related  to  Aloff--Wallach seven--manifolds. In Section \ref{section_de2} we shall introduce the generalised
Wilczynski invariants, and show
(Theorem \ref{theo_lift}) that  only cubics give rise to a complete analytic family on a contact complex three--fold $\PP(T\CP^2)$.

The connection between the algebraic geometry of cuspidal cubics, the differential geometry of their parameter space, and the 7th order differential equations
 forms a part of more general theory \cite{DG10}. Any seven--dimensional 
family of rational curves 
which lifts to a complete family of non--singular contact curves in contact complex three--folds gives rise to a $G_2$ structure. This structure in general has torsion which can be expressed in terms of contact invariants of the associated 7th order ODE characterising the curves. 
\vskip5pt
{\bf Acknowledgements.}  Some of this work was carried over
when both authors visited the Institute of Mathematics of the Polish Academy of Sciences in Warsaw and the Erwin Schr\"odinger Institute in Vienna. We thank both institutions   
for hospitality and financial support.
The second author (MD) is grateful to Burt Totaro for 
useful discussions.
\section{Group actions on symmetric polynomials}
Most calculations in the paper rely on an explicit description of the 
actions of the projective linear group on spaces
of homogeneous polynomials, and in this Section 
we shall summarise the basic facts
and notation here. Let $E=\C^m$ be an $m$--dimensional complex
vector space, and let $\mbox{Sym}^k (E)$ be the vector space
of complex homogeneous polynomials of degree $k$ in $m$ variables.
The two cases of interest will be the space
of binary sextic where $(k=6, m=2)$, and the space of ternary
cubics where $(k=3, m=3)$. 

The general element
of $\mbox{Sym}^k (E)$ is of the form
\[
P(Z)=\sum_{\alpha_i=1}^m P_{\alpha_1\alpha_2\dots\alpha_k} 
Z^{\alpha_1}Z^{\alpha_2}\dots Z^{\alpha_k},
\]
where the symmetric tensor  $P_{\alpha_1\alpha_2\dots\alpha_k}$
consists of the coefficients of the polynomial, and 
$Z^\alpha=[Z^1, Z^2, \dots, Z^k]$ are homogeneous coordinates
on $\PP(E)$. The linear action of $GL(E)$ on $E$ is given by ordinary matrix multiplication $Z\rightarrow\hat{Z}$, where $Z^\alpha={N^\alpha}_\beta \hat{Z}^\beta$ for
$N\in GL(E)$. This induces a nonlinear action of 
$GL(E)$ on $\mbox{Sym}^k (E)$ given by
$P(Z)=\hat{P}(\hat{Z})$, so that
\be
\label{non_linear_action}
\hat{P}_{\alpha_1\alpha_2\dots\alpha_k}= 
\sum_{\beta_i=1}^m {N^{\beta_1}}_{\alpha_1} 
{N^{\beta_2}}_{\alpha_2}\dots  {N^{\beta_k}}_{\alpha_k} 
{P}_{\beta_1\beta_2\dots\beta_k}.
\ee
This action preserves the homogeneity of the polynomials,
so it induces a nonlinear projective group action
of $PGL(E)$ on $\PP(E)$. Before going any further we shall restrict ourselves to the two cases of interest
\subsection{Binary sextics and classical invariants}
Let $E=\C^2$ and
let $\mathcal{V}_6=\mbox{Sym}^6(\C^2)$ be the seven dimensional space
of binary sextics of the form (\ref{sextic111}) with $n=6$. 
Let
\[
V(s, t)=
v_0t^6+6v_1t^5s+15v_2t^4s^2+20v_3t^3s^3+15v_4t^2s^4+6v_5ts^5+v_6s^6
\in \mathcal{V}_6.
\]
\begin{defi}
An invariant of a binary sextic
under the $GL(2, \C)$ action  (\ref{non_linear_action})
is a function $I=I(v_0, \dots, v_{6})$ 
such that
\[
I( \hat{v}_0, \dots, \hat{v}_{6})
=(\alpha\delta-\beta\gamma)^w I( v_0, \dots, v_{6} ), \quad\mbox{where}
\quad N=
\left( \begin{array}{cc}
\alpha & \beta\\
\gamma & \delta
\end{array} \right).
\]
The number $w$ is called the weight of the invariant.
\end{defi} 
One of the classical results of the invariant theory is that all invariants
arise from the {transvectants} (see e.g. \cite{Grace_Young}, \cite{olver2} ).
There are five invariants for binary sextics of degrees 2,  4, 6, 10 and 15 respectively connected by a syzygy of degree 30. The one we will be concerned with is the quadratic invariant
\be
\label{I2}
I_2(V):=<V, V>_6=v_0 v_6-6v_1 v_5+15v_2 v_4-10(v_3)^2.
\ee
Invariants of several binary forms arise in an analogous way, and
we shall make use of an invariant of three binary sextics, which should be thought of as a scalar part in the Clebsh--Gordan decomposition
of $\mathcal{V}_6\otimes \mathcal{V}_6\otimes\mathcal{V}_6$.
This is given by
\be
\label{trv3}
I_3(U, V, W):=<<U, V>_3, W>_6.
\ee
The invariant $I_2$ defines a symmetric quadratic form. The invariant $I_3$ is anti--symmetric
in any pair of vectors.
\subsection{Ternary cubics and their orbits}
Let $E=\C^3$. We shall consider the space
of irreducible ternary cubics which give rise to plane cubic curves in $\CP^2$
of the form
\[
\sum_{\alpha,\beta,\gamma=1}^3 P_{\alpha\beta\gamma} Z^\alpha Z^\beta Z^\gamma=0.
\]
There are ten coefficients $P_{\alpha\beta\gamma}$ but the overall scale
is unimportant, so the space of such cubics is
$\CP^9$. The group action (\ref{non_linear_action}) preserves
 the homogeneity, so it
descends to the projective action of $PGL(3, \C)$ on $\CP^9$.
There are three types of orbits (see e. g. \cite{harris}) which
we shall present in inhomogeneous coordinates
$x=Z^1/Z^3, y=Z^2/Z^3$.
\begin{enumerate}
\item Smooth cubic $y^2=x(x-1)(x-c)$.
\item Nodal cubic $y^2=x^3-x^2$. 
\item Cuspidal cubic $y^2=x^3$. 
\end{enumerate}
\begin{center}
\includegraphics[width=4cm,height=5cm,angle=270]{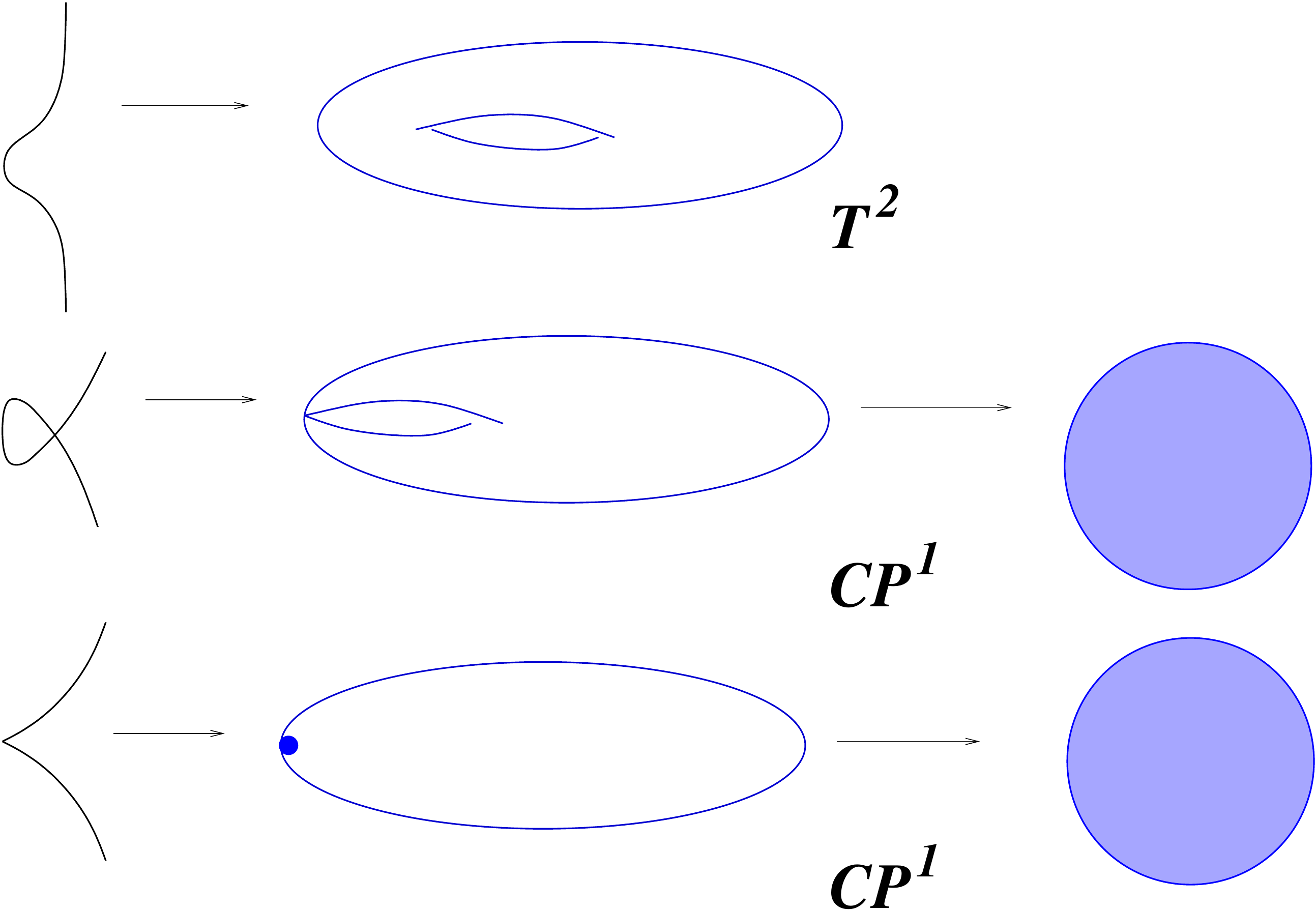}
\begin{center}
{{\bf Figure 1.} {\em Three types of orbits of irreducible cubics.}}
\end{center}
\end{center}
The first orbit corresponds to a smooth cubic which is a curve
of genus one. Two smooth cubics are projectively equivalent if
their $j$--invariants $j=(c^2-c+1)^3/(c^2(c-1)^2)$
coincide. The remaining two orbits correspond to singular rational cubics. There
is the eight dimensional orbit of the nodal cubic $(Z^2)^2Z^3-(Z^1)^3+(Z^1)^2Z^3=0$
which contains all nodal cubics, and finally there is the seven--dimensional orbit $M$ of the cuspidal cubic $Z^3(Z^2)^2-(Z^1)^3=0$.
The one--dimensional stabiliser of this cubic is given by the projective transformations with $N=\mbox{diag}(a, a^4, a^{-5})$, where $a\in\C^*$.
Thus the space of cuspidal cubics is a homogeneous manifold $M=SL(3, \C)/\C^*$.
\section{Cuspidal cubics, their moduli, and the $GL(2)$ structure}
\label{section_gl2}
\vskip3pt
In this section we shall use the twistor correspondence
to construct a $GL(2)$ structure on the space
of cuspidal cubics.
\subsection{$GL(2)$ structures}
\begin{defi}
The $GL(2)$ structure on an $(n+1)$--dimensional manifold
$M$ is an isomorphism
\be
\label{isom}
TM\cong\mbox{Sym}^n(\spp),
\ee
where $\spp$ is a rank two symplectic vector bundle over $M$.
\end{defi}
The $GL(2)$ structures were originally called the paraconformal
structures in \cite{DT}.
If $M$ is a complex manifold we talk about $GL(2, \C)$ structures
and $\spp$ is a complex vector bundle. 
The tangent vector fields to $M$ are identified by (\ref{isom})
with homogeneous polynomials of degree $n$ in two variables. There is also a unique, up to scale, symplectic structure on the fibres $\C^2$ of $\spp$.
The group action (\ref{non_linear_action}) with $m=2$ and $k=n$ gives rise
to an irreducible $(n+1)$ representation of $GL(2, \C)$,  and thus to the embedding 
of $GL(2, \C)$ inside $GL(n+1, \C)$. The image
of $SL(2, \C)\subset GL(2, \C)$ is contained in $Sp(n+1, \C)$ if $n$ is odd, or in $SO(n+1, \C)$ if $n$ is even. 
We shall consider the case of even $n$, where the representation space $\mathcal{V}_n$ is odd--dimensional. Then we have two--real sections of the $GL(2, \C)$ 
structures on real $(n+1)$--dimensional manifold $M$ which correspond to two real forms of $GL(2, \C)$.
\begin{enumerate}
\item The $GL(2, \R)$ structures give an identification
of tangent vectors with real homogeneous polynomials of degree $n$.
The fibers of $\spp$ are real vector spaces $\R^2$.
\item The $SO(3, \R)\times\R^*$ structures 
(or locally equivalent $U(2)$ structures)
identify the tangent vectors to $M$ with harmonic homogeneous polynomials
in three variables. This has its roots in 
the isomorphism $\C^3=\mbox{Sym}^2(\C^2)$ between complex vectors
in $\C^3$ and symmetric $2$ by $2$ matrices with complex coefficients.
The $(n+1)$--dimensional space $\mbox{Sym}^n(\C^2)$ 
is then identified with a subspace of $\mbox{Sym}^{n/2}(\C^3)$ which consist
of harmonic ternary forms, i. e. those forms 
$\sum_{\alpha,\dots, \gamma=1}^{n/2}P_{\alpha \beta\dots \gamma }Z^\alpha Z^\beta\dots Z^\gamma$ which satisfy $\sum_{\alpha, \beta=1}^{n/2}
\delta^{\alpha \beta} P_{\alpha\beta\dots \gamma}=0$.
\end{enumerate}
In practice
the isomorphism (\ref{isom}) is specified by a homogeneous polynomial
$S$ of degree $n$ with values in $\Lambda^1(M)$. Given
$S\in \Lambda^1(M)\otimes {\mathcal V}_n$, the homogeneous polynomial
corresponding under (\ref{isom}) to a vector field $V\in TM$ is the contraction
$V\hook S$. 
\subsection{Cuspidal cubics and their deformations}
A cuspidal cubic $L\subset \CP^2$ is a singular rational curve with
self-intersection number $9$ as  two general cuspidal cubics 
intersect in exactly nine points (albeit not in the general positions).
The arithmetic genus is constant in algebraic
families, so this is the same as the genus of a smooth cubic curve.
The arithmetic genus of a curve with a cusp is 1 plus the genus
of a resolution of singularities.) The Riemann-Roch theorem for 
singular curves yields
  \[
h^0(L, N(L)) - h^1(L,N(L)) = deg(N(L))-g(L)+1
                             = 9-1+1
                             = 9.
\]
In the case of cuspidal cubics, $h^1(L,N(L))=0$ and so $h^0(L,N(L))=9$.
Here $h^0(L,N(L))$ is equal to the dimension of the space
of all deformations of $L$ as a curve in $\CP^2$. Indeed, the space
of all cubic curves in $\CP^2$ is isomorphic to $\CP^9$. Thus 
$H^0(L,N(L))$
describes all deformations of $L$ as a curve in $\CP^2$, not just
those as a cuspidal curve.

We want to consider the deformations of $y^2-x^3=0$ as a {\em rational 
cuspidal} curve, and not allow the perturbations of cuspidal cubics to smooth curves. Thus we shall compute the Zariski tangent space to $M$ at a point in $M$. This will make use of the rational parametrisation of $L$ 
and lead to the 
the first order deformations\footnote{An alternative construction based on 
contact resolution of the cusp will be presented in Theorem \ref{theo_lift}.}.

 Let $L_m\in \mathcal{Z}=\CP^2$ be a cuspidal cubic corresponding to $m\in M$. Consider a neighboring curve in $\mathcal{Z}$ 
corresponding to a point $m+\delta m$ in $M$. In 
the proof of Proposition \ref{prop_GL} where we shall show 
that two nearby cuspidal cubics $L_m$ and $L_{m+\delta m}$ intersect 
at six points away from the triple intersection at the cusp.
This will follow from the  form of the normal vector field to $L$. 
\begin{center}
\includegraphics[width=6cm,height=8cm,angle=270]{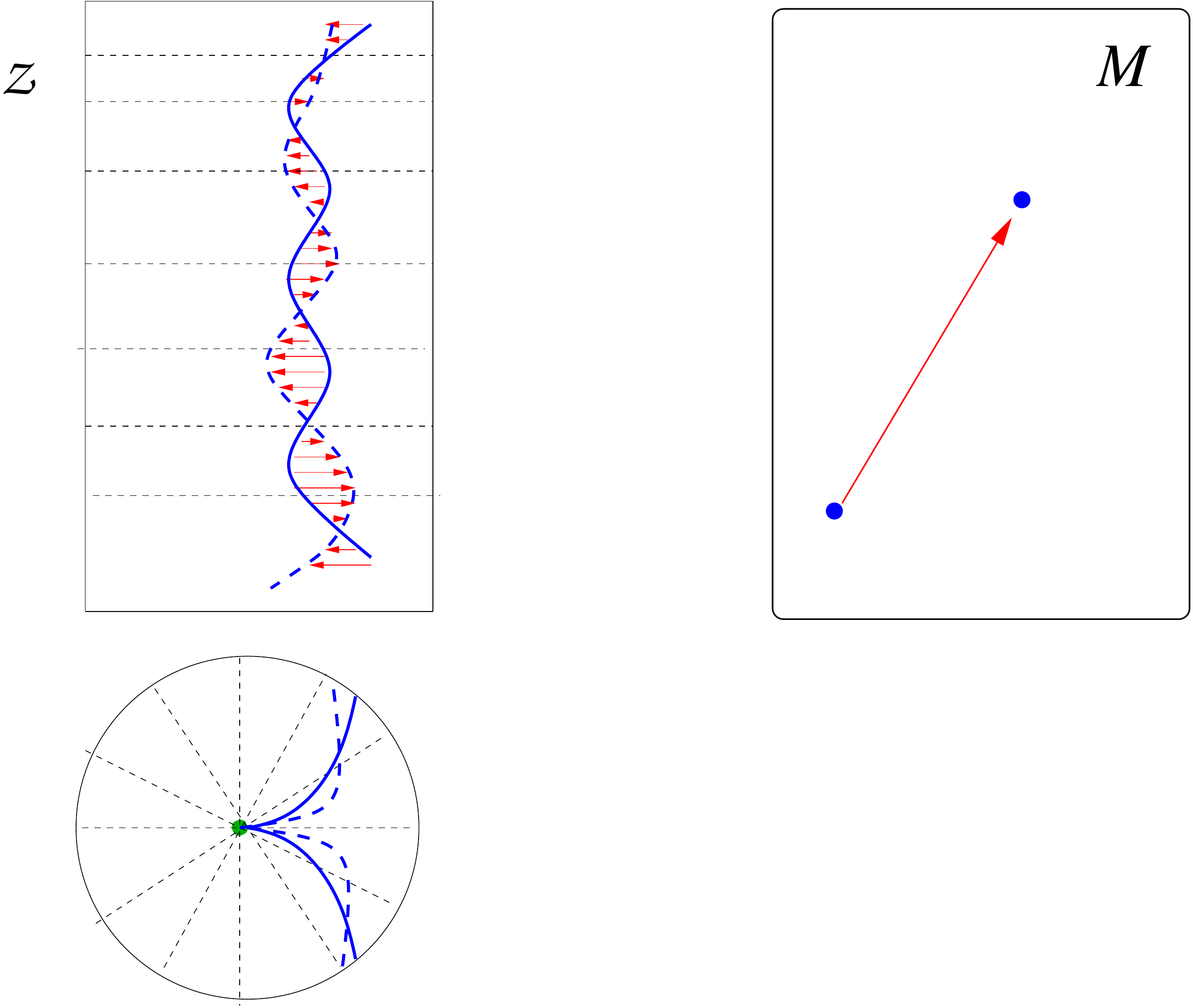}
\begin{center}
{{\bf Figure 2.} {\em Normal vector field to $L_m\subset\mathcal{Z}$ 
and tangent vector
at $m\in M$.}}
\end{center}
\end{center}
Thus a vector in $M$ connecting two neighboring points corresponds to a six--order homogeneous polynomial (defined by its six roots) in two variables which gives rise to an isomorphism
\be
\label{first_iso}
T_mM\cong\mbox{Sym}^6(\C^2).
\ee
To make contact with the notation in formula (\ref{non_linear_action})
we let $\hat{Z}^\alpha=[\hat{Z}^1, \hat{Z}^2, \hat{Z}^3]$ be homogeneous
coordinates on $\CP^2$.
\begin{prop}
\label{prop_GL}
The seven--dimensional space of cuspidal cubics 
$M=SL(3, \C)/\C^*$ admits a $GL(2)$ structure (\ref{first_iso}). If
\be
\label{ho_cubic}
\sum_{\alpha, \beta, \cdots, \phi=1}^3 P_{\alpha\beta\gamma} {N^{\alpha}}_\delta {N^{\beta}}_\epsilon
{N^{\gamma}}_\phi \;\hat{Z}^{\delta}\hat{Z}^{\epsilon} \hat{Z}^{\phi}=0,
\ee
where 
$
P_{111}=-1, \,P_{223}=P_{232}=P_{322}=1/3
$
(and other components of $P$ vanish) is the general element of $M$ then
the $GL(2)$ structure is given by $S\in \Lambda^1(M)\otimes \mbox{Sym}^6(\C^2)$
\be
\label{sextic3} 
S(s, t)=2{\sigma^2}_3\;s^6-3{\sigma^{1}}_3\;ts^5 +2{\sigma^2}_1\;t^2s^4+({\sigma^3}_3+2{\sigma^2}_2
-3{\sigma^1}_1)\;t^3s^3-3{\sigma^1}_2\;t^4s^2+{\sigma^3}_1\; t^5s
+{\sigma^3}_2\;t^6,
\ee
where $\sigma=N^{-1}dN$ is the Maurer--Cartan one--form on $SL(3, \C)$ with values
in the Lie algebra of traceless $3\times 3$ matrices.
\end{prop}
{\bf Proof.}
Consider a rational parametrisation of (\ref{cuspidal})
\[
x=t^2, \quad y=t^3,
\]
or, in homogeneous coordinates $Z^{\alpha}=T^{\alpha}(t)=[t^2, t^3, 1]$.
To establish he isomorphism (\ref{first_iso})
we shall construct a binary sextic in $(s, t)$ with values
in $T^*M$.  The $PSL(3)$ orbit of (\ref{cuspidal}) is 
seven--dimensional, and is parametrised by the components of the
matrix $N\in SL(3, \C)$. We identify two such matrices if they differ
by a multiplication by the stabiliser $\mbox{diag}(a, a^4, a^{-5})$.

Formula (\ref{non_linear_action}) with $d=3$ thus implies that the homogeneous form of the general cuspidal cubics is
(\ref{ho_cubic}), as this is the general orbit of $y^2-x^3=0$.
To construct the normal vector field to this family, differentiate (\ref{ho_cubic})
with respect to the moduli parameters ${N^{\alpha}}_{\beta}$, and 
substitute the rational parametrisation
\[
\hat{Z}^{\alpha}(t)=\sum_{\beta=1}^3{(N^{-1})^{\alpha}}_\beta T^{\beta}(t), \quad \mbox{where}\quad
T^{\alpha}=(t^2, t^3, 1).
\]
In general, if the family of rational curves
$f(x, y; m)=0$ parametrised by $m\in M$ admits a rational parametrisation
$x=x(t, m), y=y(t, m)$, then the polynomial is given by \cite{DG10}
\[
\sum_{k=1}^{\mbox{dim}\;M}\frac{\p f}{\p m^k}|_{\{x=x(t, m), y=y(t, m)\}}\; dm^k.
\]
In our case this gives a polynomial of degree nine in $t$
\be
\label{pol_sex}
\sum_{\alpha,\beta,\gamma,\delta=1}^3 P_{\alpha\beta\gamma}{\sigma^{\gamma}}_\delta T^{\alpha}T^{\beta}T^{\delta},
\ee
where $\sigma=N^{-1}dN$ is the Maurer--Cartan one--form on $SL(3, \C)$ with values
in the Lie algebra of traceless $3\times 3$ matrices. Pulling out the overall scalar factor of $t^3$ and introducing the homogeneous coordinates
$(s, t)$ in place of an affine coordinate $t$
yields the $T^*M$--valued sextic polynomial (\ref{sextic_isom}) given by 
(\ref{sextic3}). Its roots depend on coordinates on $M$.
\koniec
In particular the proof above shows  that the normal vector field to a cuspidal cubics
vanishes to the third order at the cusp $t=0$, and at six smooth points on the cubics in general positions.

\section{Construction of the $G_2$ structure}
\label{section_main}
In this section we shall present the proof of Theorem
\ref{main_theo} and show that the space of cuspidal cubics $M$ is
equipped with a conformal $G_2$ structure 
(a good reference to $G_2$ structures is \cite{salamon})  . In our discussion of $GL(2, \C)$ structures we have noted the existence of the embedding of $SL(2, \C)$ in $SO(7, \C)$. The representation theoretic argument of \cite{morozov}, or the explicit construction in \cite{DG10} shows that this leads to an intermediate embedding
$SL(2, \C)\subset {G_2}^\C \subset SO(7, C)$ which we shall now 
explore\footnote{Dynkin has shown (see e.g. \cite{morozov}) that
for general $n$ there no proper Lie subgroup  $G$ of 
$Sp(n+1, \C)$ or $SO(n+1, \C)$ such that $SL(2, \C)\subset G$.
The exception is $n=6$, where $G={G_2}^{\C}$. Thus in all dimensions apart
from seven the $GL(2)$ structure does not induce any additional $G$
structure on $M$ apart from a conformal structure if $n$ is even, 
or a symplectic structure if $n$ is odd.} .

Let $V,U,W\in TM$. The $GL(2)$ structure allows the identification
of vector fields with binary sextics, and therefore the invariants
(\ref{I2}) and (\ref{trv3}) give rise to a non--degenerate symmetric quadratic form
and a skew--symmetric three--form  on $M$ given by
\be
\label{formula_new}
g(U, V)=<U\hook S, V\hook S>_6, \quad
\phi(U, V, W)=<<U\hook S, V\hook S>_3, W\hook S>_3.
\ee
For a general $GL(2)$ structure  the isomorphism (\ref{first_iso}) is specified by a $T^*M$ valued
sextic polynomial
\be
\label{sextic_isom}
S(s, t)=a^0t^6+6a^1t^5s+15a^2t^4s^2+20a^3t^3s^3+15a^4t^2s^4+
6a^5ts^5+a^6s^6,
\ee
for linearly independent one--forms $a^0, \dots, a^6$ on $M$.
Given $S(s, t)$, a sextic polynomial corresponding to a vector
$V\in TM$ is given by the contraction $V\hook S(s, t)$. The transvectant formulae for the invariants (\ref{I2}) and (\ref{trv3}) , and formula
(\ref{formula_new}) imply that the quadratic form $g$ and the three--form
$\phi$ are given by
\be
\label{conf_inv}
g= a^0\odot a^6-6a^1\odot a^5+15a^2\odot a^4-10(a^3)^2,
\ee
and 
\be 
\label{form_inv}
\phi={\sqrt\frac{5}{2}}\Big(3\,(a^1\wedge a^2\wedge a^6+a^0\wedge a^4\wedge a^5)+
a^3\wedge(a^0\wedge a^6+6\;a^1\wedge a^5-15\, a^2\wedge a^4)\Big).
\ee
The quadratic invariant $g$ induces a conformal structure on $M$ and
$\phi$  (the overall multiple $\sqrt{5/2}$ has been chosen for later 
convenience) endows $M$ with a three--form compatible with $g$ in a sense
that 
\[
(V\hook \phi)\wedge (V\hook \phi)\wedge \phi=0 \qquad\mbox{iff}\qquad g(V, V)=0.
\]
The invariants (\ref{I2}) and (\ref{trv3}) have weight six and nine respectively,
and thus $g$ gives rise to a conformal structure, but not a metric.
Changing a metric in the conformal structure has to be complemented
by changing the three--form according to
\[
g\longrightarrow  \Omega^6 g, \quad
 \phi\longrightarrow 
\Omega^9 \phi,
\]
where $\Omega$ is a non--vanishing function on $M$. Thus the structure group of $TM$ reduces to the complexification of conformal $G_2$ \cite{DG10}.
We shall now find this structure explicitly, and demonstrate that it is co--calibrated.
 
{\bf Proof of Theorem \ref{main_theo}.} Using the form of the $GL(2)$ structure
given by (\ref{sextic3}) together with formulae (\ref{conf_inv}) and 
(\ref{form_inv}) gives rise to a conformal structure represented by the metric
\be
\label{metric_2}
g=2{\sigma^3}_2\odot{\sigma^2}_3+\frac{1}{2}{\sigma^3}_1\odot{\sigma^1}_3
-\frac{2}{5}{\sigma^1}_2\odot{\sigma^2}_1-\frac{1}{40}(4{\sigma^1}_1-
{\sigma^2}_2)^2.
\ee
This holomorphic conformal structure on  $SL(3, \C)/\C^*$ admits
three real forms,\textit{\textit{}} all leading to co-calibrated $G_2$ structures.
\begin{itemize}
\item If the components of $\sigma$ are all real, then $g$ is a 
metric of signature $(3, 4)$ on the non--compact manifold
$SL(3, \R)/\R^*$.
\item If $\sigma$ is anti--hermitian then $g$ is a metric of signature
$(4, 3)$ on a compact seven--manifold $M=SU(3)/U(1)$, where $U(1)$ is the group 
of diagonal matrices  
\be
\label{u1stab}
\left( \begin{array}{ccc}
e^{i\theta} & 0&0\\
0 & e^{4i\theta}&0\\
0 & 0& e^{-5i\theta}
\end{array} \right), \quad \theta\in\R.
\ee
Thus
$M$ is the Aloff--Wallach space $N(1, 4)$. The homogeneous co--calibrated 
$G_2$ structures of Riemannian signature on some Aloff--Wallach spaces have appeared before in \cite{swann,Gibbons,ilka,R10,Baum}.
\item If the diagonal components of $\sigma$ are imaginary and
the reality conditions
\be
\label{reality}
 {\sigma^2}_1=- \overline{{\sigma^1}_2}, \quad
{\sigma^3}_1=\overline{{\sigma^1}_3},\quad
{\sigma^3}_2=\overline{{\sigma^2}_3}
\ee
hold, then the metric (\ref{metric_2}) has Riemannian signature. The relations
(\ref{reality}) imply that $\sigma$ takes values in $\mathfrak{su}(2, 1)$. Thus
we obtain a homogeneous Riemannian conformal class on the non--compact 
seven--manifold $M=SU(2, 1)/U(1)$, where $U(1)$ is given by (\ref{u1stab}).
\end{itemize}
We shall now work out the details of the $G_2$ structure (\ref{metric_2})
associated with the Riemannian reality conditions (\ref{reality}), and show that it is co--calibrated with the conformal factor equal to 1.
Let 
\[
e_8=\left( \begin{array}{ccc}
i & 0&0\\
0 & 4i&0\\
0 & 0& -5i
\end{array} \right)
\]
span the one dimensional Lie algebra of the $U(1)$ stabiliser 
(\ref{u1stab}). We choose the following basis for the invariant complement of
$e_8$ (the various square roots multiples are chosen such that the dual  
one--forms in the resulting metric have length one)
\begin{eqnarray*}
e_3&=&\frac{\sqrt{10}}{2} \left( \begin{array}{ccc}
0 & -1&0\\
1 & 0&0\\
0 & 0& 0
\end{array} \right), \quad
e_7= \frac{\sqrt{10}}{2}\left( \begin{array}{ccc}
0 & i&0\\
i & 0&0\\
0 & 0& 0
\end{array} \right),\quad
e_2=\sqrt{2}\left( \begin{array}{ccc}
0 & 0&1\\
0 & 0&0\\
1 & 0& 0
\end{array} \right),\\
e_6&=&\sqrt{2}\left( \begin{array}{ccc}
0 & 0&-i\\
0 & 0&0\\
i & 0& 0
\end{array} \right),\quad
e_1=\frac{1}{\sqrt{2}}\left( \begin{array}{ccc}
0 & 0&0\\
0 & 0&1\\
0 & 1& 0
\end{array} \right),
\quad
e_5=\frac{1}{\sqrt{2}}\left( \begin{array}{ccc}
0 & 0&0\\
0 & 0&i\\
0 & -i& 0
\end{array} \right),\\
e_4&=&\frac{\sqrt{10}}{7}\left( \begin{array}{ccc}
-3i & 0&0\\
0 & 2i&0\\
0 & 0& i
\end{array}\right).
\end{eqnarray*}
The Maurer--Cartan one--form on $SU(2, 1)$ is
\[
\sigma=N^{-1}dN = \left( \begin{array}{ccc}
{\sigma^{1}}_1 & {\sigma^{1}}_2  & {\sigma^{1}}_3\\
{\sigma^{2}}_1  & {\sigma^{2}}_2  &{\sigma^{2}}_3 \\
 {\sigma^{3}}_1 & {\sigma^{3}}_2 & -{\sigma^{1}}_1 -{\sigma^{2}}_2
\end{array} \right)=
\sum_{k=1}^8 e_k\otimes\theta^k,
\]
where $\theta^k$ are the left--invariant one--forms on the group.
Thus 
\[
{\sigma^1}_2=\frac{\sqrt{10}}{2}(-\theta^3+i\theta^7), 
\quad {\sigma^1}_3=\sqrt{2}(\theta^2-i\theta^6),
\quad {\sigma^2}_3=\frac{1}{\sqrt{2}}(\theta^1+i\theta^5),\quad 
{\sigma^2}_2-4{\sigma^1}_1=2i\sqrt{10}\;\theta^7
\]
together with the reality conditions (\ref{reality}). The metric
(\ref{metric_2}) and the $G_2$ three--form (\ref{form_inv}) become
\begin{eqnarray}
\label{formula1}
g&=&(\theta^1)^2+ (\theta^2)^2+(\theta^3)^2+(\theta^4)^2+(\theta^5)^2
+(\theta^6)^2+(\theta^7)^2,\\
\phi&=&\theta^{123}+ \theta^{145}+ \theta^{167}
+ \theta^{246}- \theta^{257} -\theta^{347} - \theta^{356}\nonumber,
\end{eqnarray}
where $\theta^{jkl}=\theta^j\wedge \theta^k \wedge \theta^l$.
The 
relations 
\[
d\sigma+\sigma\wedge\sigma=0, \qquad [e_j, e_k]=\sum_{l=1}^3 c_{jkl}\;e_l
\]
give
\begin{eqnarray*}
d\theta^1&=&\sqrt{10}\theta^2\wedge\theta^3+\frac{\sqrt{10}}{7}\theta^4\wedge\theta^5
-9\theta^5\wedge\theta^8+\sqrt{10}\theta^6\wedge\theta^7,\\
d\theta^2&=&-\frac{\sqrt{10}}{4}\theta^1\wedge\theta^3 +\frac{4\sqrt{10}}{7}\theta^4\wedge\theta^6
-\frac{\sqrt{10}}{4}\theta^5\wedge\theta^7+6\theta^6\wedge\theta^8,\\
d\theta^3&=&-\frac{\sqrt{10}}{5}\theta^1\wedge\theta^2+ \frac{5\sqrt{10}}{7}\theta^4\wedge\theta^7
+\frac{\sqrt{10}}{5}\theta^5\wedge\theta^6-3\theta^7\wedge\theta^8,\\
d\theta^4&=&\frac{\sqrt{10}}{20}\theta^1\wedge\theta^5 +\frac{4\sqrt{10}}{5}
\theta^2\wedge\theta^6
-\frac{5\sqrt{10}}{4}\theta^3\wedge\theta^7,\\
d\theta^5&=&\frac{\sqrt{10}}{7}\theta^1\wedge\theta^4+9\theta^1\wedge\theta^8+
\sqrt{10}\theta^2\wedge\theta^7+\sqrt{10}\theta^3\wedge\theta^6,\\
d\theta^6&=&-\frac{\sqrt{10}}{4} \theta^1\wedge\theta^7 +\frac{4\sqrt{10}}{7}\theta^2\wedge\theta^4
-6\theta^2\wedge\theta^8 -\frac{\sqrt{10}}{4}\theta^3\wedge\theta^5,\\
d\theta^7&=&-\frac{\sqrt{10}}{5}\theta^1\wedge\theta^6+
\frac{\sqrt{10}}{5}\theta^2\wedge\theta^5+
\frac{5\sqrt{10}}{7}\theta^3\wedge\theta^4 +3\theta^3\wedge\theta^8,\\
d\theta^8&=&\frac{3}{14}\theta^1\wedge\theta^5 -\frac{4}{7}\theta^2\wedge\theta^6-\frac{5}{14}\theta^3\wedge\theta^7
\end{eqnarray*}
and finally (\ref{cocalibrated}).
Thus the $G_2$ structure is co-calibrated.
\koniec
The real form $SU(2, 1)/U(1)$ which we have explored in the proof corresponds
to the second `real form` of the $GL(2, \C)$ structure (see the discussion
of the real forms  in Section \ref{section_gl2}). Thus it is an example of the 
$SO(3)$ structure in seven dimensions \cite{fried_new}.

A similar construction applied to nodal cubics would instead lead to a symplectic
structure on the eight dimensional $PSL(3)$ orbit. To make contact
with conformal geometry one needs to blow up a point 
in $\CP^2$ to get a seven--dimensional family $M$. This will admit a 
non--homogeneous  $G_2$ structure. 
\section{ODEs for cuspidal curves}
\label{section_de1}
All rational curves with seven--dimensional orbits are projectively 
equivalent to the cuspidal curves
\be
\label{new_cuspics}
y^p-x^q=0
\ee
with $(p, q)$ integers. This fact follows from a more general result established in \cite{Doubrov3}. In this  section
we shall present a direct proof based on the projective curvature. 
We shall also  characterise  the family 
of cuspidal cubics and their higher degree generalisations
(\ref{new_cuspics}) as integral 
curves of 7th order ODEs.

 In Wilczy\'nski's approach to projective differential geometry \cite{Wilczynski_book}
each curve $\C\rightarrow \CP^2$ (or $\R\rightarrow \RP^2$) corresponds to
a unique 
third order homogeneous linear ODE
\be
\label{Wilczynski}
\frac{d^3 Y}{dx^3}+3p_1(x)\frac{d^2 Y}{dx^2}+3p_2(x)\frac{d Y}{d x}+p_3(x)Y=0
\ee
such that given a curve 
$x\rightarrow [y_1(x), y_2(x), y_3(x)]$, the functions ${\bf y}=[y_1(x), y_2(x), y_3(x)]$ 
span the solution space of (\ref{Wilczynski}). To find this ODE,
substitute each $y_i(x)$ into (\ref{Wilczynski}) and solve 
the resulting system of linear algebraic equations for each of the smooth functions $p_i$s.
Linear transformations of the basis ${\bf y}$ correspond to  projective 
transformations of the curve. These transformations do not change the ODE
(\ref{Wilczynski}). The combinations of the coefficients which only 
depend on the ratios of the solutions , i.e. are unchanged by
transformations ${\bf y}\rightarrow \gamma(x){\bf{y}}$
are the semi--invariants 
\be
\label{P2P3}
P_2=p_2-(p_1)^2-(p_1)_x, \quad P_3=p_3-3p_1p_2+2(p_1)^3-(p_1)_{xx},
\ee
where the subscripts stand for partial derivatives. The lowest order relative
projective invariant is given by
\[
\Theta_3(x)=P_3-\frac{3}{2}(P_2)_x.
\]
The cubic differential
$
\Theta_3(x) dx^3
$
is invariant under an overall scaling of homogeneous coordinates and 
reparametrisation of the curve,
\[
(x, {\bf y})\longrightarrow (\xi(x), \gamma(x){\bf y}).
\]
Invariants of the ODE (\ref{Wilczynski}) under this 
class of transformations are also projective invariants of the curve.
Using these  transformations one can set two out of the three functions 
$p_i$ to zero.
\begin{itemize}
\item The Laguerre--Forsyth canonical form is achieved by setting
$p_1=p_2=0$ in which case  $\Theta_3=p_3$.
Thus if $\Theta_3=0$, the solution space is
${\bf y}=[1, x, x^2]$ and the curve is a conic.
\item Consider a curve $y=y(x)$, so that ${\bf y}=[1, x, y(x)]$. This gives
$p_2=p_3=0$ and $p_1=-y_{xxx}/(3y_{xx})$. Then
\be
\label{conic_ode}
\Theta_3=\frac{9(y^{(2)})^2y^{(5)}-45y^{(2)}y^{(3)}y^{(4)}+40(y^{(3)})^3}
{{(y^{(2)})}^3}.
\ee
This, as we have just shown,  vanishes for conics 
which gives a characterisation of the five dimensional space of plane conics by a 5th order ODE
originally due to Halphen \cite{Halphen}. \end{itemize}
We will say that $\Theta_r(x)$ is a relative invariant of weight $r$ if
$\Theta_r(x) dx^r$ is an invariant. Wilczynski shows that given $\Theta_r$,
the quantity
\be
\label{Thetar}
\Theta_{2r+2}=2r\Theta_r{(\Theta_r)}_{xx}-(2r+1)({(\Theta_r)}_x)^2-3r^2P_2(\Theta_r)^2
\ee
is a relative invariant of weight $(2r+2)$. Thus $\Theta_3$ gives rise
to $\Theta_8$, and we can define the projective curvature to be an absolute invariant
\be
\label{curvature}
\kappa=\frac{(\Theta_8)^3}{(\Theta_3)^8}.
\ee
This is the lowest order absolute 
projective invariant\footnote{The fact that there are no invariants
of order lower than seven can also be seen by direct counting. The prolongations
of the eight generators of $PSL(3)$ from the $(x, y)$ plane to the
6th jet $J^6$ are vector fields which are independent almost everywhere, and
thus span $TJ^6$ at almost every point. Therefore the only functions
of $(x, y, y', \dots, y^{(6)})$ constant along the flows generated by the lifts
are constant identically.}. If a parametrisation with $p_2=p_3=0$ is 
chosen, then the expression for $\kappa$ depends on $y$ and its first seven derivatives. 

\begin{prop}[Wilczynski \cite{Wilczynski_book}, Sylvester \cite{sylvester}]
\label{cubic_ode_prop}
The cuspidal cubics have constant projective curvature, and are 
characterised by the 7th order ODE
\be
\label{ccode}
\kappa(y, y', \cdots, y^{(7)})=\frac{3^9 7^3}{2^4 5^2}  .
\ee
\end{prop}
{\bf Proof.}
This can be seen directly parametrising the cuspidal cubic by 
$[1,x, x^{3/2}]$ so that  equation (\ref{Wilczynski}) becomes
\[
\frac{d^3 Y}{d x^3}+\frac{1}{2x}\frac{d^2 Y}{d x^2}=0,
\]
and $\kappa$ can be found directly by substituting $p_1=1/(6x)$ in 
the formulae above. 
\koniec

It is easy to find all other curves with constant projective curvature. A curve $[1, x, x^\gamma]$, where $\gamma\neq 0, 1$ 
is characterised by
\[
\kappa=3^9\frac{(1+\gamma^2-\gamma)^3}{(\gamma-2)^2(2\gamma-1)^2(\gamma+1)^2}.
\]
The case 
\[
\kappa=3^9/2^2
\]
which corresponds to $\gamma=0$ or  $\gamma=1$ has to be considered separately,
as for these two values the solutions to Wilczy\'nski's ODE  (\ref{Wilczynski})
are not independent. We verify that this special case corresponds to a
curve $[1,x, \ln(x)]$. Therefore we have established
\begin{prop}
\label{flat_curves}
All curves of constant projective curvature
are projectively equivalent to 
\[
y=x^\gamma,\;\gamma\neq 0, 1, -1, 2, 1/2 \quad  \mbox{or}\quad y=\ln{x}.
\]
\end{prop}
This is in agreement with \cite{Doubrov3}, where the same class of curves arose
as the most general homogeneous curves in $\CP^2$.
All algebraic curves in this class are projectively equivalent to 
the rational cuspidal curves (\ref{new_cuspics}).
The stabiliser of (\ref{new_cuspics}) is the one--dimensional 
group of matrices
\[
\left( \begin{array}{ccc}
a^{q-2p}&0&0\\
0&a^{p-2q} &0\\
0&0 &a^{p+q}
\end{array} \right), \quad a\in \C^*.
\]
There are three real forms of the space of orbits, 
as before. One of these is the Aloff-Wallach space $N(k, l)=SU(3)/U(1)$, where the $U(1)$ 
subgroup consists of matrices of the form \cite{AW}
\[
\left(\begin{array}{ccc}
e^{il\theta}&0&0\\
0&e^{ik\theta}&0\\
0&0&e^{-i(k+l)\theta}
\end{array} \right),\quad \theta\in \R 
\]
and $k, l$ are integers such that
\[
p=\frac{2l+k}{3}, q=-\frac{l+2k}{3}.
\]
The space $N(1, 1)$ is special from this perspective. It corresponds to
curves of the form $xy=1$. These curves   form a five--dimensional orbit 
$SL(3)/SL(2)$ which is the space  of all conic sections.

\section{Generalised Wilczynski invariants for the 
constant projective curvature}
\label{section_de2}
\subsection{Classical Wilczynski invariants}
Equation~\eqref{curvature} with constant $\kappa$ is a 7th order ODE whose solutions are curves with constant projective curvature $\kappa$. Therefore, for any given $\kappa$ we have a seven--dimensional space of such curves and each of the curves has a one-dimensional stabiliser.  All these solution spaces can be identified with homogeneous spaces $PSL(3,\R)/H$, where $H$ is one of the one-dimensional subgroups in $PSL(3,\R)$ described above.

In ~\cite{DG10} it has been demonstrated that a solution space to a 7th order ODE admits a  $G_2$-structure if all its {\it generalised Wilczynski invariants} vanish identically. To introduce these invariants recall that given  an arbitrary linear differential equation
\begin{equation}\label{lin_eq}
Y^{(n)}+\binom{n}{1}p_1(x)Y^{(n-1)}+\binom{n}{2}p_2(x)Y^{(n-2)}+\dots+p_n(x)Y(x)=0,
\end{equation}
with real smooth or complex holomorphic coefficients $p_i(x)$, $i=1,\dots,n$ the classical Wilczynski invariants are constructed as follows: First, we can always bring this equation to the 
so-called semi-canonical form
\begin{equation}\label{semi-form}
Y^{(n)}+\binom{n}{2}P_2(x)Y^{(n-2)}+\binom{n}{3}P_3(x)Y^{(n-3)}+\dots+P_n(x)Y(x)=0.
\end{equation}
This is achieved by
\begin{equation}\label{semi-transf}
Y\mapsto \lambda Y,\quad \text{where }\lambda = \exp(-\int p_1(x)dx). 
\end{equation}
It is easy to check that the new coefficients $P_i(x)$, $i=2,\dots,n$ are polynomial expressions in terms of $p_i(x)$, $i=1,\dots,n$ and their derivatives. They also do not depend on the integration constant in ~\eqref{semi-transf}. For example $P_2$ and $P_3$ are
given by (\ref{P2P3}).
Next, the semi-canonical form~\eqref{semi-form} can be brought to the 
 Laguerre--Forsyth canonical form: 
\begin{equation}\label{canonical-form}
Y^{(n)}+\binom{n}{3}q_3(x)Y^{(n-3)}+\dots+\binom{n}{n-1}q_{n-1}(x)Y'(x)+q_n(x)Y(x)=0,
\end{equation}
by means of the following change of variables:
\begin{equation}\label{canonical-transf}
(x,Y) \mapsto ( \xi(x), (\xi')^{(n-1)/2}Y ),
\end{equation}
where $\xi$ is satisfies the third order ODE reduced to the Riccati equation
\begin{equation}\label{riccati}
\eta'-1/2\eta^2 = \frac{6}{n+1}P_2,\quad \eta = \xi''/\xi'.
\end{equation}
In general, the coefficients $q_i(x)$, $i=3,\dots,n$ of the canonical form~\eqref{canonical-form} depend on the the choice of the solution $\xi(x)$ of defined by~\eqref{riccati}. 
\begin{theo}[Wilczynski \cite{Wilczynski_book}] The expressions:
\begin{equation}\label{wil_inv}
\Theta_r =  \frac{1}{2}\sum_{s=0}^{r-3}(-1)^s\frac{(r-2)!r!(2r-s-2)!}{(r-s-1)!(r-s)!(2r-3)!s!}q_{r-s}^{(s)}
\end{equation}
are relative invariants of the linear differential 
equation (\ref{lin_eq}) with respect to the transformations
$(x,y)\mapsto (\xi(x),\lambda(x)y )$, i. e.  $\Theta_r \mapsto (\xi')^r\Theta_r$.
\end{theo}
The expressions~\eqref{wil_inv} were introduced by Wilczynski in~\cite{Wilczynski_book} and are called the
linear (relative) invariants of the equation~\eqref{lin_eq}. Wilczynski also shows that all other invariants (relative and absolute) can be defined from them via differentiation and algebraic operations. In particular, equation~\eqref{lin_eq} can be transformed to the trivial equation $Y^{(n)}=0$ if and only if $\Theta_r=0$ for all $r=3,\dots,n$. 

We note  that the invariants $\Theta_r$ do not depend on the choice of $\eta$ in~\eqref{riccati} and can be expressed explicitly in terms of coefficients $p_i(x)$, $i=1,\dots,n$, of the initial equation~\eqref{lin_eq}. Do do it in practice, bring (\ref{lin_eq})  to the form (\ref{semi-form}) and calculate the coefficients $q_i$
using (\ref{canonical-transf}). Now $q_1$ vanishes identically, but $q_2$ does not unless
the Riccati equation holds. We nevertheless formally compute the expressions~\eqref{wil_inv} with substitutions $\xi''=\xi'\eta$ and $\eta'=1/2\eta+\frac{6}{n+1}P_2$. The coefficients in \eqref{wil_inv} are chosen in such a way that the  resulting expression will no longer depend on~$\eta$. For example, the explicit expression for $\Theta_3$ does not depend on the order $n$ and has the form:
\[
\Theta_3 = p_3 - 3p_1p_2+2p_1^3+3p_1p_1'-\frac{3}{2}p_2'+\frac{1}{2}(p_1')^2.
\]
For $n=3$ we arrive at the invariant $\Theta_3$ defined in the previous section. 

Other approaches to define Wilczynski invariants~\eqref{wil_inv} are presented in works~\cite{ovsienko,chalkley}. In particular, R.~Chalkley provides an algorithm for computing Wilczynski invariants in a way that avoids the Laguerre--Forsyth canonical form. He also gives alternative proofs of the above results.

We note that Wilczynski invariants can also be computed in the cases when the coefficients of the initial equation have isolated singularities. In particular, they are well-defined meromorphic functions if the linear
equation is defined on complex domain with all coefficients being meromorphic functions on this domain. 

\subsection{Generalized Wilczynski invariants} 
The generalised Wilczynski 
invariants~\cite{Doubrov1,Doubrov2}
of an arbitrary non-linear ODE
\be
\label{ODE_n}
y^{(n)}=F(x,y,y',\dots,y^{(n-1)})
\ee
are defined as classical Wilczynski invariants of its linearisation. Analytically, they are computed by substituting $-\binom{n}{r}^{-1} D_x^k\big(\frac{\partial F}{\partial y^{(n-r)}}\big)$ in place of $p_r^{(k)}(x)$ in the classical Wilczynski invariants. Here by $D_x$ we denote the operator of total derivative:
\[
D_x = \frac{\partial}{\partial x} + y' \frac{\partial }{\partial y} + \dots + y^{(n-1)}\frac{\partial}{\partial y^{(n-2)}} +  F\frac{\partial}{\partial y^{(n-1)}}.
\]
We denote by ${\bf\Theta}_i$ the generalised Wilczynski invariant we obtain from $\Theta_i$ by this formal substitution.

\begin{theo}[\cite{Doubrov1}]
Generalised Wilczynski invariants are (relative) contact invariants of non-linear ordinary differential equations.
\end{theo}

Generalised Wilczynski invariants were defined and studied in~\cite{Doubrov1, Doubrov2}. They are closely related to W\"unschmann conditions defined 
in~\cite{DT} and further explored in \cite{GN}. 
For $n=7$ these W\"unschmann conditions are explicitly computed in~\cite{DG10}. They consist of five expressions $W_1,\dots, W_5$, which should vanish identically to guarantee the existence of a natural $GL(2, \R)$-structure 
(\ref{isom}) on the solution space\footnote{Not all $GL(2)$ structures 
come from ODEs. Those which do  have been partially characterised in \cite{Kry}.}
$M$ of the equation, such that normal vector to a hypersurface given by fixing $(x, y)$ in the general solution to (\ref{ODE_n}) corresponds to 
a sextic polynomial with a 
root of multiplicity $6$. For $n=7$ these W\"unschmann conditions are related to the generalised Wilczynski invariants as follows:
\[
W_i = a_i \left(\mathbf\Theta_{i+2} + \sum_{j=3}^{i+1} b_i^j\mathbf\Theta_j\right),
\]
where $a_i$ are constants and $b_i^j$ are linear differential operators which are polynomials in the total derivative operator $D_x$ of order at most $4$. Eg., $a_1=-3430$, $a_2=-240100$ and $b_2^3 = \tfrac{2}{5}D_x - \tfrac{12}{35}\tfrac{\partial F}{\partial y^{(6)}}$ so that:
\begin{align*}
W_1 &= -3430\; {\bf\Theta}_3,\\
W_2 &= -240100\left( {\bf \Theta}_4 + \tfrac{2}{5}D_x({\bf\Theta}_3) - \tfrac{12}{35}\tfrac{\partial F}{\partial y^{(6)}}{\bf \Theta}_3\right).
\end{align*}
In particular, it immediately follows that vanishing of ${\bf\Theta}_i$, $i=3,\dots,7$ is equivalent to vanishing of $W_i$, $i=1,\dots,5$.

\begin{lemma} 
\label{lemma_ku}
Let $(\Theta_8)^3 = \kappa (\Theta_3)^8$ be the
 differential equation (\ref{curvature}) defining all curves on projective plane with the constant projective curvature $\kappa\neq 0$. All generalized Wilczynski invariants
vanish for this equation if and only if $\kappa = \frac{3^9 7^3}{2^4 5^2}$. It corresponds exactly to the family of all cuspidal cubics.
\end{lemma}
{\bf Proof.}
The explicit expression for $\Theta_3$ is given by equation~\eqref{conic_ode}. The explicit expression for $\Theta_8$ is given by~\eqref{Thetar} and is linear in the highest derivative $y^{(7)}$. Resolving the 
equation $(\Theta_8)^3 = \kappa (\Theta_3)^8$ with respect $y^{(7)}$ and using the explicit formulae for Wilczynski invariants~\eqref{wil_inv}, we get that the generalised Wilczynski invariants for equation~\eqref{curvature} are given by
\begin{align*}
\mathbf\Theta_3 &=\mathbf\Theta_4=\mathbf\Theta_5=\mathbf\Theta_7= 0,\\
\mathbf\Theta_6 &= -\frac{(2^4 5^2 \kappa - 3^9 7^3)}{2^23^{12} 7^4} \frac{\bigg(9(y^{(2)})^2y^{(5)}-45y^{(2)}y^{(3)}y^{(4)}+40(y^{(3)})^3\bigg)^2}{(y'')^6}.
\end{align*}
Thus, we see that $\kappa = \frac{3^9 7^3}{2^4 5^2}$ is the only value of projective curvature, for which all generalised Wilczynski invariants vanish identically.
\koniec
\subsection{Orbits of general cuspidal curves}
\begin{theo}
\label{theo_lift}
Let $\mathcal{C}_{(p, q)}=SL(3, \C)/\C^*$ be a seven--dimensional family of all plane curves projectively equivalent to the curve
\[
y^p-x^q=0,
\]
where $(p, q)$ are coprime positive integers with $\mbox{Max}(p, q)>2$. The following statements
are equivalent
\begin{enumerate}
\item
The seven--dimensional family $\mathcal{C}_{(p, q)}$ is a complete contact deformation family of non-singular Legendrian curves in $\PP(T\CP^2)$.
\item
Generalised Wilczynski invariant of the ODE (\ref{curvature}) vanish.
\item
$\mathcal{C}_{(p, q)}$ is projectively equivalent to the family of cuspidal cubics $\mathcal{C}_{(2, 3)}$.
\end{enumerate}
\end{theo}
{\bf Proof.}
Let us consider the $SL(3)$ orbits of general cuspidal curves~\eqref{new_cuspics}. The reason for
non--vanishing of $\mathbf\Theta_6$ is that the $SL(3)$ orbit is not a complete analytic
family in the sense  of \cite{kodaira} unless the curve is cubic. Repeating the steps
leading to (\ref{sextic3}), with $T^{\alpha}=(t^p, t^q, 1)$ gives
\begin{eqnarray*}
S(s, t)&=&
(q-p){\sigma^3}_2\; t^{2q}+(q-p){\sigma^3}_1\;t^{q+p}s^{q-p}
-q{\sigma^1}_2\;t^{2q-p} s^p\\
&&+((q-p){\sigma^3}_3+p{\sigma^2}_2
-q{\sigma^1}_1)\;t^qs^{q} +p{\sigma^2}_1\;
t^ps^{2q-p}-q{\sigma^1}_3\;t^{q-p}s^{q+p}+p{\sigma^2}_3\;s^{2q}.
\end{eqnarray*}
The form of this polynomial is not preserved by the rational transformation
of $t$ unless Max$(p, q)=3$.

Let also clarify this fact from the point of view of singularities of Legendrian curves. Passing to the homogeneous coordinates $[X: Y: Z]$ on $\CP^2$ and permuting the coordinates if needed, we can always assume that the equation of the curve is written as $Y^pZ^{q-p}=X^q$, where $p\le q$ and $q\ge 3$. In this case the point $[0: 0: 1]$ is always singular, so that all curves in our family will be singular as well. 

In order to apply Bryant's generalisation \cite{Bryant2} of Kodaira deformation theory we need all curves of the family $\mathcal{C}_{(p,q)}$ to be non-singular, or at least their lifts to the projectivized cotangent bundle $\PP(T\CP^2)$  to be non-singular. Elementary computation shows that the lift of the curve $Y^pZ^{q-p}=X^{q}$ to $P(T\CP^2)$ is a non-singular curve if and only if $p=1$ or $q=p+1$. Indeed,
 the parametrisation of the affine coordinates$(x, y)=(t^p, t^q)$ lifts to a rational parametrisation of the
curves $\gamma(t)=(x, y, \zeta)=(t^p, t^q, (q/p)t^{q-p})$ 
which are Legendrian with respect to the contact form $dy-\zeta dx$.
We find that $\gamma=\dot{\gamma}=0$ at $t=0$ unless $p=1$, or $q=p+1$.

Both cases are projectively equivalent via the change of coordinates $Y$ and $Z$. So, we shall treat only the first case and assume that $p=1$. Then only for $q=3$ (the cuspidal cubics case) the lifts of all curves in the family $\mathcal{C}_{(p,q)}$ constitute the complete deformation family of the curve $YZ^2=X^3$. Indeed, Lemma \ref{lemma_ku} shows that  for $q>3$ the projective curvature~\eqref{curvature} is different from the distinguished value~$\frac{3^97^3}{2^45^2}$ and the generalized Wilczynski invariant ${\bf\Theta}_6$ does not vanish for the 7th order ODE defining the family of curves projectively equivalent to $YZ^{q-1}=X^q$.
\koniec

The fact that the generalised Wilczynski invariant ${\bf\Theta}_6$ does not vanish on curves $YZ^{q-1}=X^q$ for $q > 3$ although all these curves are rational can be explained as follows. Consider 
the 7th order ODE defining the family of curves $\mathcal{C}_{(1,q)}$.
Formulae (\ref{Thetar}) with $r=3$ and~\eqref{curvature} imply 
that the coefficient of the leading term $y^{(7)}$ equal exactly to the nominator of
the classical Wilczynski invariant $\Theta_3$ (see~\eqref{conic_ode}). Thus, the equation of plane conics
\[
9(y^{(2)})^2y^{(5)}-45y^{(2)}y^{(3)}y^{(4)}+40(y^{(3)})^3=0
\]
defines at the same time the set of singular points for the equation of curves of constant projective curvature.
This is exactly the set where all generalized Wilczynski invariants have singularities as well. Now computing the 
expression  (\ref{conic_ode}) for the curve $YZ^{q-1}=X^q$, $q\ge 3$, we find that it is equal to $X^{3q-9}$ up to a non-zero constant. Thus, the lifts of curves projectively equivalent to $YZ^{q-1}=X^q$ cross the set of singular points transversally whenever $q>3$ and do not intersect this set at all if $q=3$. In other words, the generalized Wilczynski invariants have no singularities on the curves $YZ^{q-1}=X^q$ if and only if $q=3$.

\section{Outlook. ODEs for rational curves.}
We have constructed a co-calibrated $G_2$ structure on the moduli space $M$
of cuspidal cubics.  Co-calibrated $G_2$ structures play a role in
theoretical 
physics:
they give rise to  solutions of IIB supergravity for which
the only flux is the self-dual five-form \cite{Gutowski1}.
They also appear in the context of near-horizon geometries in
heterotic supergravity. In particular, the eight--dimensional
spatial cross-sections of the horizon are $U(1)$ fibrations over
a conformaly co-calibrated $G_2$ structures on a seven--manifold $M$ 
\cite{Gutowski2}. There is also a connection with $SU(3)$ structures
\cite{chiosi}.

In our work, the manifold $M$  arises as the solution space
of the 7th order ODE (\ref{ccode}). This is an example of the general 
construction of \cite{DG10} which associates $G_2$ structures with  $7$th
order ODEs with general solutions given by rational curves. 

There are few known ODE with that property - they are partially
characterised by the vanishing of the Wilczynski invariants
\cite{Doubrov1, DT, Doubrov2, GN}. They also correspond to projective differential invariants \cite{olver}.
\begin{itemize} 
\item The ODE 
\[
25y^{(7)} (y^{(4)})^2-105 {y^{(6)}y^{(5)}}{y^{(4)}}+{84}{(y^{(5)})^3}=0
\]
describes rational sextics with two cusps and admits eight--dimensional group of symmetries (this group is different than $PSL(3)$). The corresponding
$G_2$ structure is closed
\[
d\phi=0, \quad d*\phi= \tau\wedge \phi,
\]
where $\tau$ is some two--form  \cite{DG10}.

\item 
The ODE
\[
10 (y^{(3)})^3y^{(7)}- 70(y^{(3)})^2y^{(4)}y^{(6)}
-49(y^{(3)})^2(y^{(5)})^2
+280(y^{(3)})(y^{(4)})^2y^{(5)}
-175 (y^{(4)})^4=0
\]
is (together with $y'''=0)$ the unique ODE admitting ten dimensional group
of contact symmetries \cite{noth}. The general solution is given by certain 
family of rational sextics \cite{DS10}.
The symmetry group is isomorphic to $Sp(2)$, and the seven--dimensional solution space $M=Sp(2)/SL(2)$ admits  a nearly--integrable $G_2$  structure 
\[
d\phi=\lambda *\phi, \qquad d*\phi=0,
\] 
where $\lambda$ is a constant. If the real form $SO(5)/SO(3)$ is chosen,
then the $G_2$ structure is Riemannian \cite{Bryant1}.
\end{itemize}
There are also some lower order examples. 
\begin{itemize}
\item The 5th order ODE for conics 
\[
9(y^{(2)})^2y^{(5)}-45y^{(2)}y^{(3)}y^{(4)}+40(y^{(3)})^3=0
\]
is equivalent to the vanishing of  $\Theta_3$  given 
by (\ref{conic_ode}). This ODE goes back at least to Halphen \cite{Halphen}. 
The resulting solution space $SL(3)/SL(2)$ admits a homogeneous metric which can be read off from the general conic in a way analogous to our 
construction of (\ref{metric_2}). 
\item Considering all conics 
passing through two points
$[1, 0, 0], [0, 1, 0]$ in $\CP^2$ gives a three--dimensional family
\[
y=\frac{ax+b}{cx+d}.
\]
This arises from the Schwartzian ODE
\[
2y^{(3)}y^{(1)}-3(y^{(2)})^2=0.
\]
\end{itemize}
The problem of classifying ODEs whose all solutions are rational curve
remains open.

\end{document}